\documentclass{article}
\usepackage{amsmath,amsfonts,latexsym,amssymb,amsthm}
\usepackage{url}

\begin{document}

\title{The Diophantine equation $x^4\pm y^4=iz^2$ in Gaussian integers}
\author{Filip Najman}
\date{}
\maketitle
\begin{abstract}
In this note we find all the solutions of the Diophantine equation $x^4\pm y^4=iz^2$ using elliptic curves over $\mathbb Q(i)$. Also, using the same method we give a new proof of Hilbert's result that the equation $x^4\pm y^4=z^2$ has only trivial solutions in Gaussian integers.
\end{abstract}

\section{Introduction.}
The Diophantine equation $x^4\pm y^4=z^2$, where $x,y,$ and $z$ are integers, was studied by Fermat, who proved that there exist no nontrivial solutions. Fermat proved this using the \emph{infinite descent} method, proving that if a solution can be found, then there exists a smaller solution (see for example \cite[Proposition 6.5.3]{coh}). This was the first particular case of Fermat's Last Theorem proven (the theorem was completely proven by Wiles in \cite{wil}).

The same Diophantine equation, but now with $x,y,$ and $z$ being Gaussian integers, i.e., elements of $\mathbb Z[i]$, was later examined by Hilbert (see \cite[Theorem 169]{hil}). Once again, it was proven that there exist no nontrivial solutions. Other authors also examined similar problems. In \cite{sz} equations of the form $ax^4+by^4=cz^2$ in Gaussian integers with only trivial solutions were studied. In \cite{cr} a different proof than Hilbert's is given, using descent, that $x^4+y^4=z^4$ has only trivial solutions in Gaussian integers. The equations $x^4\pm y^4=z^2$ over some quadratic fields were also considered in \cite{sz2} and \cite{sz3}, again proving that there exist no nontrivial solutions. Some applications of Diophantine equations of this type can be found in \cite{xu1} and \cite{xu2}. 

In this short note, we will examine the Diophantine equation 
$$x^4\pm y^4=iz^2$$
in Gaussian integers and find all the solutions of this equation. Also, we will give a new proof of Hilbert's results. Our strategy will differ from the one used by Hilbert and will be based on elliptic curves. Elliptic curves have also been used in \cite{lam} to prove that the Diophantine equation $x^3+y^3=z^3$ has only trivial solutions in Gaussian integers, but in a somewhat different way than in this note. We will use elliptic curves over quadratic fields that have nontrivial torsion, while in \cite{lam}, an elliptic curve with trivial torsion over the rationals was examined.

For an elliptic curve $E$ over a number field $K$, it is well known, by the Mordell-Weil theorem, that the set $E(K)$ of points on the curve $E$ with coordinates in the field $K$ is a finitely generated abelian group. The group $E(K)$ is isomorphic to $T\oplus\mathbb Z^r$, where $r$, which is called the \emph{rank}, is a nonnegative integer and $T$, which is called the \emph{torsion subgroup}, is the group of all elements of finite order. Thus, there are finitely many points on an elliptic curve over a field if and only if it has rank 0.
   
We will be interested in the case when $K=\mathbb Q(i)$. We will work only with elliptic curves with rational coefficients, and by a recent result of the author (see \cite{fn}), if an elliptic curve has rational coefficients, then the torsion of the elliptic curve over $\mathbb Q(i)$ is either cyclic of order $m$, where $1 \leq m \leq 10$ or $m=12$, of the form $\mathbb Z_2 \oplus \mathbb Z_{2m}$, where $1 \leq m \leq 4$, or $\mathbb Z_4 \oplus \mathbb Z_4$.

Throughout this note, the following extension of the Lutz-Nagell Theorem is used to compute torsion groups of elliptic curves. 
\newtheorem*{tm100}{Theorem (Extended Lutz-Nagell Theorem)}
\begin{tm100}
Let $E$ be an elliptic curve in the form $E: y^2=x^3+Ax+B$ with $A,B\in \mathbb Z[i]$. If a point $(x,y)\in E(\mathbb Q(i))$ has finite order, then:
\begin{enumerate}
\item $x,y\in \mathbb Z[i].$
\item Either $y=0$ or $y^2|4A^3+27B^2.$
\end{enumerate} 
\end{tm100}
The main step in the proof of the Lutz-Nagell Theorem (for curves over $\mathbb Q$) is to show that all the torsion points have integer coordinates. This is done by showing that no prime can divide the denominators of the coordinates of the torsion points. The proof of the Lutz-Nagell Theorem can easily be extended to elliptic curves over $\mathbb Q(i)$. Details of the proof can be found in \cite[Chapter 3]{th}. An implementation in Maple can be found in \cite[Appendix A]{th}.

Note that every elliptic curve can be put in the form $E: y^2=x^3+Ax+B$ (this is the \emph{short Weierstrass form}) over any field of characteristic zero, and thus in particular over $\mathbb Q(i)$. The Extended Lutz-Nagell Theorem enables us to easily get a finite list of possible candidates for the torsion points, and then check which ones are actually torsion points. Although when all the torsion points are found, one could easily compute the group structure of the torsion subgroup using addition laws, all we need is the list of all the torsion points.\\

\vspace{0.1cm}

\section{The Diophantine equation $x^4\pm y^4=iz^2$.}
\theoremstyle{definition}
\newtheorem{tm5}{Definition}
\begin{tm5}
We call a solution $(x,y,z)$ of the Diophantine equation 
$$x^4\pm y^4=cz^2,$$
for some given $c\in \mathbb C$, \emph{trivial} if $xyz=0$.
\end{tm5}
We are now ready to prove our main result. 
\theoremstyle{plain}
\newtheorem{tm2}{Theorem}
\begin{tm2}
\begin{itemize}
\item[(i)] The equation $x^4-y^4=iz^2$ has only trivial solutions in Gaussian integers.
\item[(ii)] The only nontrivial solutions satisfying $\gcd(x,y,z)=1$ in Gaussian integers of the equation $x^4+y^4=iz^2$ are $(x,y,z)$, where $x,y\in \{\pm i, \pm 1\},\ z=\pm i(1+i)$.
\end{itemize}
\end{tm2}
\emph{Proof:}\\
(i) Suppose $(x,y,z)$ is a nontrivial solution. Dividing the equation by $y^4$ and making the variable change $s= x/y,\ t= z/y^2$, we obtain the equation
$s^4-1=it^2$, where $s,t\in \mathbb Q(i)$. We can rewrite this equation as 
\begin{equation}
r=s^2,
\end{equation}
\begin{equation}
r^2-1=it^2.
\end{equation}
Multiplying these equations we obtain $i(st)^2=r^3-r$. Again, making the variable change $a=st,\ b=-ir$ and dividing by $i$, we obtain an equation defining an elliptic curve $$E:a^2=b^3+b.$$ Using the program \cite{sim}, written in PARI, we compute that the rank of this curve is 0. It is easy to compute, using the Extended Lutz-Nagell Theorem, that $E(\mathbb Q(i))_{tors}=\mathbb Z_2\oplus\mathbb Z_2$ and that $b\in\{0,\pm i\}$. It is obvious that all the possibilities lead to trivial solutions.\\
(ii) Suppose $(x,y,z)$ is a nontrivial solution satisfying $\gcd(x,y,z)=1$. Dividing the equation by $y^4$ and making the variable change $s= x/y,\ t=z/y^2$, we obtain the equation
$s^4+1=it^2$, where $s,t\in \mathbb Q(i)$. We can rewrite this equation as 
\begin{equation}
r=s^2,
\label{eq1}
\end{equation}
\begin{equation}
r^2+1=it^2.
\end{equation}
Multiplying these equations we obtain $i(st)^2=r^3+r$. Again, making the variable change $a=st,\ b=-ir$ and dividing by $i$, we obtain an equation defining an elliptic curve $$E:a^2=b^3-b.$$ Using the program \cite{sim}, we compute that the rank of this curve is 0. Using the Extended Lutz-Nagell Theorem we compute that $E(\mathbb Q(i))_{tors}=\mathbb Z_2\oplus\mathbb Z_4$ and that $b\in\{0,\pm i,\pm 1\}$. Obviously $b=0$ leads to a trivial solution. It is easy to see that $b=\pm 1$ leads to $r=\pm i$ and this is impossible, since $r$ has to be a square by (\ref{eq1}). This leaves us the possibility $b=\pm i$. Since we can suppose that $x$ and $y$ are coprime, this case leads us to the solutions stated in the theorem. \qed\\

\section{A new proof of Hilbert's results.}
We now give a new proof of Hilbert's result, which is very similar to Theorem 1.
\newtheorem{tm}[tm2]{Theorem}
\begin{tm}
The equation $x^4\pm y^4=z^2$ has only trivial solutions in Gaussian integers.
\end{tm}
\emph{Proof:}\\
Suppose $(x,y,z)$ is a nontrivial solution. Dividing the equation by $y^4$ and making the variable change $s= x/y,\ t=z/y^2$, we obtain the equation
$s^4\pm 1=t^2$, where $s,t\in \mathbb Q(i)$. We can rewrite this equation as 
\begin{equation}
r=s^2,
\end{equation}
\begin{equation}
r^2\pm 1=t^2,
\end{equation}
and by multiplying these two equations, and making the variable change $a=st$, we get the two elliptic curves
$$a^2=r^3\pm r.$$
As in the proof of Theorem 1, both elliptic curves have rank 0 and it is easy to check that all the torsion points on both curves lead to trivial solutions.\qed\\
\theoremstyle{remark}
\newtheorem*{rem}{Remark}
\begin{rem}
Note that from the proofs of Theorems 1 and 2 it follows that the mentioned solutions are actually the only solutions over $\mathbb Q(i)$, not just $\mathbb Z[i]$.
\end{rem}

\paragraph{Acknowledgments.}  The author would like to thank the referees for many helpful suggestions. The author was supported by the Ministry of Science, Education and Sports, Republic of Croatia, Grant 037-0372781-2821.

\bigskip

\noindent\textit{Department of Mathematics, University of Zagreb, Bijeni\v cka cesta 30, 10000 Zagreb, Croatia\\
fnajman@math.hr}
\end{document}